\documentclass[12pt]{article}

\usepackage{hyperref,yfonts,pifont,mathrsfs,latexsym,palatino,amsthm,amssymb,amsmath,amsfonts,graphicx,color,tikz}

\allowdisplaybreaks

\newtheorem{theorem}{Theorem}[]

\newtheorem{remark}[theorem]{Remark}

\newtheorem{proposition}[theorem]{Proposition}

\makeatletter
\newcommand{\leqnomode}{\tagsleft@true\let\veqno\@@leqno}
\newcommand{\reqnomode}{\tagsleft@false\let\veqno\@@eqno}
\makeatother

\begin{document}

\section*{Cantor's Non-Equinumerosity Theorems, Inductively}
\vspace{-.30cm}

\begin{center}
\textbf{Saeed Salehi}\footnote{%
Harish and Bina Shah School of AI and Computer Science, Plaksha University, IT City Road, Sector 101A, Mohali, Punjab 140306, India. \textbf{\texttt{root@SaeedSalehi.ir}}}
\end{center}

\begin{quote}{\small
\textbf{Abstract.}
In the first, pre-college level part, we present a proof by mathematical induction for Cantor’s theorem on the uncountability of the infinite binary strings and the real numbers. In the second, undergraduate-level part, we prove Cantor’s powerset theorem by using transfinite induction. There, we will need the axiom of choice and the concept of ordinals. Some of these proofs by (mathematical and transfinite) induction seem to be new; comparisons will be made with older results.

{\bf AMS 2020 Subject Classification}:
  03E10 
   $\cdot$
 03E25

   {\bf Keywords} ({\sl in alphabetic order}):
  Axiom of Choice $\cdot$
      Binary Expansions $\cdot$
  Binary Strings $\cdot$
    Mathematical Induction  $\cdot$
    Ordinal Numbers $\cdot$
    Power-Set $\cdot$
 Real Numbers $\cdot$
   Transfinite Induction $\cdot$
   Uncountability $\cdot$
  Well-Order
 }\end{quote}

\section{Introduction}
\label{sec:intro}

\noindent
{\sc  Georg Cantor} (1845---1918) is known as the founding father of modern set theory; {\sl Cantor’s paradise} refers to his naive set theory. As a matter of historical fact, Cantor was an expert in mathematical analysis, and his inventions were in response to his needs in studying the behavior of some functions on real numbers. Some theorems are named after him, such as (1) the uncountability of the set of all infinite binary strings, (2) the  uncountability of the set of all real numbers, and (3) the
non-equinumerosity of a set with its powerset;  we will study these theorems in this paper. Cantor proved them by what is now called “Cantor’s diagonal argument” (the uncountability of the reals was proved by him twice before the discovery of the diagonal argument). We will give alternative proofs by induction for these three theorems. The uncountability of the infinite binary string and the reals can be proved by using mathematical induction. Mathematical induction can also prove the uncountability of the powerset of natural numbers. To prove the general powerset theorem (that no set can surjectively cover its powerset), we will need a more general induction, so-called {\sl transfinite induction}, which applies to well-ordered sets, and this requires the axiom of choice, one of whose equivalent statements is that every set can be well-ordered.

\section{\textsc{Mathematical Induction}}
\subsection{Binary Sequences}
Cantor's celebrated diagonal argument constructs a new infinite binary sequence $c$ for a given infinite sequence $\{s_i\}_{i=1}^{\infty}$ of infinite binary sequences, such that $c$ is different from each $s_i$, for $i=1,2,\cdots$. The sequence $c$ differs from the given list $s_i$ at the $i$th place, as it is defined by $c|_i={\tt 1}-s_i|_i$, where $s|_i$ denotes the $i$th element of the sequence $s$; so, we have $c|_i={\tt 0}$ when $s_i|_i={\tt 1}$ and $c|_i={\tt 1}$ otherwise (when $s_i|_i={\tt 0}$).  Here, we present a novel proof that inductively constructs a new sequence out of the given list. For an infinite sequence $s$ and a positive integer $n$, let $s|_{\leqslant n}$ denote the sequence $\langle s|_i\rangle_{i\leqslant n}=\langle s|_1,s|_2,\cdots,s|_n\rangle$, which is the initial segment of $s$ with length $n$. Define the binary sequence $\sigma=\langle\sigma_n\rangle_{n>0}$ by (strong) induction on $n>0$ as follows:   $$\sigma_n=\begin{cases}
                                         {\tt 0}, & \mbox{if } \langle\sigma_i\rangle_{i<n}\ast\langle{\tt 0}\rangle\neq s_n|_{\leqslant n}; \\
                                         {\tt 1}, & \mbox{if } \langle\sigma_i\rangle_{i<n}\ast\langle{\tt 0}\rangle= s_n|_{\leqslant n};
                                       \end{cases}$$
where $\ast$ denotes the concatenation of sequences. Notice that $\langle\sigma_i\rangle_{i<1}$ is empty, and so
$$\sigma_1=\begin{cases}
                                         {\tt 0}, & \mbox{if } s_1|_{\leqslant 1}\neq \langle{\tt 0}\rangle; \\
                                         {\tt 1}, & \mbox{if } s_1|_{\leqslant 1}=
                                         \langle{\tt 0}\rangle;
                                       \end{cases}
 \;\; ={\tt 1}-s_1|_1.$$

\begin{theorem}[Cantor, 1890]
\,

\noindent
  The set of infinite binary sequences is not countable.
\end{theorem}

\begin{proof}

\noindent
  For a given infinite sequence $\{s_i\}_{i=1}^{\infty}$ of infinite binary sequences, the binary sequence $\sigma=\langle\sigma_n\rangle_{n>0}$ constructed above is different from each $s_n$ since $\sigma|_{\leqslant n}=\langle
  \sigma_i\rangle_{i\leqslant n}\neq
  s_n|_{\leqslant n}$ holds for every $n>0$.
  \end{proof}

In other words, $\sigma$ differs from $s_n$ not necessarily at the $n$th place, but the initial segments of the two sequences (of length $n$) differ.

\subsection{Real Numbers}
Here, we apply this inductive argument to show the uncountability of the reals.
Fix a base $\mathfrak{b}>1$.  Every real number $r$ in the interval $(0,1]$ has a unique infinite expansion in base $\mathfrak{b}$, which is of the form $r=\sum_{i=1}^{\infty}d_i\mathfrak{b}^{-i}$, where   $d_i\in\{0,1,\cdots,\mathfrak{b}-1\}$ is not an ultimately $0$ sequence. The ultimately $0$, or finite, expansions are not unique; e.g.,   $(0.1000\cdots)_{\mathfrak{b}}
=(0.0\mathfrak{b}'\mathfrak{b}'\mathfrak{b}'\cdots)_{\mathfrak{b}}$
 holds for $\mathfrak{b}'=\mathfrak{b}-1$, since $\mathfrak{b}^{-1}=\mathfrak{b}'\sum_{i=2}^{\infty}\mathfrak{b}^{-i}$.
The diagonal argument builds a real number $r$ in $(0,1]$ for a given infinite sequence $\{r_i\}_{i>0}$ of reals in $(0,1]$, such that $r$ is different from each $r_i$ ($i>0$). To make sure that we do not produce an ultimately $0$ sequence, we take the base $\mathfrak{b}$ to be greater than $2$.
The binary base, $\mathfrak{b}=2$, causes some problems that require more care; see \textbf{\cite[{\rm p.~54}]{F76}}.
The now-classical (anti-)diagonal construction goes as follows.
Write the infinite expansion of each $r_n$ in base $\mathfrak{b}$ as $r_n=\sum_{i=1}^{\infty}d_{n,i}\mathfrak{b}^{-i}$, where $d_{n,i}$'s belong to $\{0,1,\cdots,\mathfrak{b}-1\}$ and are not ultimately $0$. Let $c=\sum_{i=1}^{\infty}c_i\mathfrak{b}^{-i}$, where the sequence $c_i$ is defined as follows: $$c_i=\begin{cases}
                1, & \mbox{if } d_{i,i}\neq 1; \\
                2, & \mbox{if } d_{i,i} = 1.
              \end{cases}$$
 The real number $c$ belongs to $(0,1]$ and can be shown to be different from each $r_n$.

\begin{theorem}[Cantor, 1874]
\,

\noindent
The real interval $(0,1]$ is uncountable.
\end{theorem}

\begin{proof}

\noindent
 Let $\{r_i\}_{i>0}$ be a given infinite sequence  of reals in $(0,1]$. Fix a base $\mathfrak{b}>2$,
and write the infinite expansion of  $r_i$ in base $\mathfrak{b}$ as $r_i=\sum_{j=1}^{\infty}d_{i,j}\mathfrak{b}^{-j}$.
For each $n,i>0$, let $r_i|_{\leqslant n}$ be the number $\sum_{j=1}^{n}d_{i,j}\mathfrak{b}^{-j}$, the restriction of the infinite expansion of $r_i$ in base $\mathfrak{b}$ to the first $n$ summands.
Define the sequence $\sigma_n$ by (strong) induction on $n>0$ as follows: $$\sigma_n=\begin{cases}
                                  1, & \mbox{if } (\sum_{i<n}\sigma_i\mathfrak{b}^{-i}+
                                  \mathfrak{b}^{-n})
                                  \neq r_n|_{\leqslant n}; \\
                                  2, & \mbox{if } (\sum_{i<n}\sigma_i\mathfrak{b}^{-i}+
                                  \mathfrak{b}^{-n}) = r_n|_{\leqslant n}.
                                \end{cases}$$
Put $\sigma=\sum_{i=1}^{\infty}\sigma_i\mathfrak{b}^{-i}$. This number is in $(0,1]$ but does not belong to $\{r_i\}_{i>0}$ since for each $n>0$, we have $\sigma|_{\leqslant n}\neq r_n|_{\leqslant n}$ by the definition of $\sigma$.
  \end{proof}

To see that this new inductive argument is really different from Cantor's diagonal argument,  let us consider   the following sequence of real numbers in $(0,1]$ over  a fixed base $\mathfrak{b}>2$:

\medskip

\qquad $r_1=0.1222222\cdots=0.1\overline{2}$,

\qquad $r_2=0.1111111\cdots=0.\overline{1}$,

\qquad $r_3=0.2222222\cdots=0.\overline{2}$,

\qquad $r_4=0.2121212\cdots=0.\overline{21}$,

\qquad $r_5=0.2111111\cdots=0.2\overline{1}$,

\qquad $r_6=0.1212121\cdots=0.\overline{12}$,

 \qquad \qquad \qquad $\vdots$

\medskip

\noindent
The (anti-)diagonal number is $c=0.221221\cdots$, while the inductive argument produces $\sigma=0.211121\cdots$.

\subsubsection{Binary Expansions}
Now, we consider the case of $\mathfrak{b}=2$.
Fix a   sequence $\{r_i\}_{i>0}$  of reals in $(0,1]$ and their infinite binary expansions $r_i=\sum_{j=1}^{\infty}b_{i,j}2^{-j}$, where $b_{i,j}\in\{0,1\}$ are not ultimately $0$.
Neither the diagonal argument nor our inductive one can work for  the binary base, since the new expansion could be ultimately $0$. To overcome this obstacle,   an ingenious argument was presented in \textbf{\cite{H10}}, which can be rephrased  as follows.

Put $h=\sum_{n=1}^{\infty}h_n2^{-n}$, where
$$h_{2n-1}h_{2n}=\begin{cases}
                   01, & \mbox{if } b_{n,2n}=0; \\
                   10, & \mbox{if } b_{n,2n}=1.
                 \end{cases}$$
Note that the binary expansion of $h$ cannot be ultimately $0$ (or $1$).
Here are two more alternative arguments.
\begin{itemize}
  \item[({\sc i})] Let for every $n>0$, $$s_{2n-1}s_{2n}=\begin{cases}
                         01, & \mbox{if } b_{n,2n-1}b_{n,2n}\neq 01; \\
                         10, & \mbox{if } b_{n,2n-1}b_{n,2n} = 01;
                       \end{cases}$$
and $s=\sum_{n=1}^{\infty}s_n2^{-n}$.
  \item[({\sc ii})] Inductively, put for every $n>0$,
$$\sigma_{2n-1}\sigma_{2n}=\begin{cases}
                         01, & \mbox{if } (\sum_{i=1}^{2n-2}\sigma_i2^{-i} + 2^{-2n}) \neq r_n|_{\leqslant 2n}; \\
                         10, & \mbox{if } (\sum_{i=1}^{2n-2}\sigma_i2^{-i} + 2^{-2n}) = r_n|_{\leqslant 2n}.
                       \end{cases}$$
The real number $\sigma=\sum_{n=1}^{\infty}\sigma_n2^{-n}$ belongs to the  interval $[1/3 , 2/3]$ and is different from each $r_n$, since $\sigma|_{\leqslant 2n}
= \sum_{i=1}^{2n}\sigma_i2^{-i}
 \neq
r_n|_{\leqslant 2n}$ holds for each $n>0$.
\end{itemize}

\section{\textsc{Transfinite Induction}}
Before going transfinite, let us try our hands on (strong) mathematical induction to prove the uncountability of the powerset of the set of all natural numbers, $\mathbb{N}$.
For a subset $S\subseteq\mathbb{N}$ and  number $n\in\mathbb{N}$, let $S_{<n}$ denote $S\cap\{0,1,\dots,n-1\}$ and $S_{\leqslant n}$ denote $S\cap\{0,1,\dots,n\}$.
For a given function $F\colon\mathbb{N}\rightarrow\mathscr{P}(\mathbb{N})$, we construct a set $X$ out of $F$'s range by induction.
To have $X\neq F(n)$ for each $n\in\mathbb{N}$, we guarantee the inequality $X_{\leqslant n}\neq F(n)_{\leqslant n}$.
At stage $0$, let $0\in X$ if and only if $0\not\in F(0)$. At stage $n>0$, put $n$ in $X$ when and only when  $X_{<n}=F(n)_{\leqslant n}$. Notice that $X_{<n}=F(n)_{\leqslant n}$ makes sense even for $n=0$, since $X_{<0}=\emptyset$ and $F(n)_{\leqslant 0}=F(0)\cap\{0\}$. Now, $X$ is not in the range of $F$, since our construction assures   $X_{\leqslant m}\neq F(m)_{\leqslant m}$ for each  $m\in\mathbb{N}$. To see this, consider two cases for a fixed $m\in\mathbb{N}$. (\texttt{\textbf{i}}) $m\in X$. Then $X_{<m} = F(m)_{\leqslant m}$ by the construction, and also, $X_{\leqslant m}\neq X_{<m}$ by the assumption $m\in X$; therefore, $X_{\leqslant m}\neq F(m)_{\leqslant m}$.  (\texttt{\textbf{i\!i}}) $m\not\in X$. Then $X_{<m}\neq F(m)_{\leqslant m}$ by the construction, and also, $X_{\leqslant m}=X_{<m}$ by the assumption $m\not\in X$; therefore, $X_{\leqslant m}\neq F(m)_{\leqslant m}$.


\subsection{The Powerset}
The inductive argument can be generalized to show the non-existence of a surjection $F\colon A\rightarrow\mathscr{P}(A)$, from a set $A$ to its powerset $\mathscr{P}(A)$, when $A$ is well-orderable.
For a subset $S\subseteq A$ of an ordered set $(A,\prec)$ and an element $a\in A$, let $S_{\prec\,a}$ denote the subset  $\{x\in S\mid x\prec a\}$ and $S_{\preccurlyeq\,a}$ denote
 $\{x\in S\mid x\preccurlyeq a\}$.

\begin{theorem}[Powerset Theorem for Well-Ordered Sets]
\,

\noindent
  If $A$ is well-orderable, then no function $F\colon A\rightarrow\mathscr{P}(A)$ can be surjective.
\end{theorem}

\begin{proof}

\noindent
Let $\prec$ be a well-ordering on $A$, and let a function $F\colon A\rightarrow\mathscr{P}(A)$ be given. Define  $B_a\subseteq A$ for $a\in A$ by (strong) $\prec$-induction as follows: $$B_a=\begin{cases}
                  \emptyset, & \mbox{if }  \bigcup_{i\prec a}B_i \neq F(a)_{\preccurlyeq\,a}; \\
                 \{a\}, & \mbox{if }  \bigcup_{i\prec a}B_i = F(a)_{\preccurlyeq\,a}.
               \end{cases}$$
Let $B=\bigcup_{a\in A}B_a$. The set $B$ is a subset of $A$ that is not in the range of $F$, since $B_{\preccurlyeq\,a}\neq F(a)_{\preccurlyeq\,a}$ holds for every $a\in A$. It suffices to notice that
$B_{\prec\,a}=\bigcup_{i\prec a}B_i\subseteq A_{\prec\,a}$ and
$B_{\preccurlyeq\,a}=\bigcup_{i\preccurlyeq a}B_i\subseteq A_{\preccurlyeq\,a}$.
\end{proof}

Let $V\sqsubseteq W$ mean that the well-ordered set $V$ is an initial segment (or equal to the whole) of the well-ordered set $W$.
The following inductive proof stems from \textbf{\cite{B97}}.

\begin{theorem}[Cantor's Powerset Theorem, \texttt{\textbf{I}}]
\,

\noindent
For a function $F\colon A\rightarrow\mathscr{P}(A)$, any $\sqsubseteq$-maximal well-orderable subset $B\subseteq A$ whose order $\prec$   satisfies the following $(\boldsymbol\ast)$ is not in the range of $F$.
$$(\boldsymbol\ast)\;\forall b\in B: F(b)=B_{\prec\,b}.$$
\end{theorem}

\begin{proof}

\noindent
Assume for the sake of a contradiction that $F(a)=B$ for some $a\in A$. Then $a\not\in B$, since otherwise, if $a\in B$, then $a\in F(a)$, and so by $\boldsymbol\ast$ we would have $a\prec a$, a contradiction. Now, $B'=B\cup\{a\}$ is bigger than $B$, and $\prec'$ that extends $\prec$ by $b\prec'a$ for every $b\in B$, well-orders $B'$, and satisfies $\boldsymbol\ast$ as well.
  This contradicts the $\sqsubseteq$-maximality of $B$.
\end{proof}

\begin{remark}{\rm
In the above proof, the condition  $\boldsymbol\ast$ could be relaxed to $$(\boldsymbol\star)\; \forall b\in B: F(b)\subseteq B_{\prec\,b}.$$

It goes without saying that the existence of such a $\sqsubseteq$-maximal subset of $A$ requires Zorn's Lemma, or its equivalent, AC. The above proof is inductive since the set $B$ can be defined as $B=\{b_{\gamma}\}_{\gamma\in\textsf{Ord}}$, where $\textsf{Ord}$ is the class of all ordinals and $b_\gamma$'s are defined by (transfinite) induction as follows: $b_\gamma$ is a member of $A$, if any, such that $F(b_\gamma) = \{b_\alpha\mid\alpha\in\gamma\}$; if there is no such member, then $b_\gamma$ is not defined for  $\gamma$ (neither for any bigger ordinal).
\hfill $\diamond$
}\end{remark}

Another inductive proof goes along the same lines by taking $(I,\prec,B)$ to be a $\sqsubseteq$-maximal well-ordered set $I$ with a mapping  $B\colon I\rightarrow\mathscr{P}(A)$ (denoted $i\mapsto B_i\subseteq A$, for each $i\in I$),  in a way that
$$(\boldsymbol\star')\; \forall j\in I\;\forall b\in B_j: F(b)\subseteq\bigcup_{i\prec j}B_i.$$
Here, $(I,\prec,B)\sqsubseteq (I',\prec',B')$ means that $(I,\prec)$ is an initial segment of, or equal to,  $(I',\prec')$ and   $B_i'=B_i$ for each $i\in I$.
For such a $\sqsubseteq$-maximal $(I,\prec,B)$, the set  $\bigcup_{i\in I}B_i$  is not in the range of $F$, since otherwise  the bigger set $I'=I\cup\{\bot\}$ (for a new $\bot\not\in I$) with $I\prec\bot$ and $B_\bot=\{a\}$ would satisfy $\boldsymbol\star'$ as well, contradicting the $\sqsubseteq$-maximality of $I$.
The main idea of the following proof is from \textbf{\cite{V13}}.

\begin{theorem}[Cantor's Powerset Theorem, \texttt{\textbf{II}}]
\,

\noindent
For a function $F\colon A\rightarrow\mathscr{P}(A)$, the set $\mathcal{B}=\bigcup_{\gamma\in\textsf{Ord}}B_\gamma$ is not in the range of $f$, where $B_\gamma$'s are defined inductively by
$$(\boldsymbol\ast')\;B_\gamma=\{x\in A\mid F(x)\subseteq\bigcup_{i\in\gamma}B_i\}.$$
\end{theorem}

\begin{proof}

\noindent
If $\mathcal{B}=F(a)$ for some $a\in A$, then $a\not\in\mathcal{B}$, since otherwise $a\in B_\mu$ for the $\in$-minimum $\mu\in{\sf Ord}$ would imply  $F(a)\subseteq\bigcup_{i\in\mu}B_i$   by $\boldsymbol\ast'$,   so $a\in F(a)$ would imply $a\in B_i$ for some $i\in\mu$, contradicting the $\in$-minimality  of $\mu$; thus, $a\not\in B_\gamma$ for each $\gamma\in{\sf Ord}$. In particular, for each ordinal $\gamma$, we have $a\not\in B_{\gamma+1}$, and so by  $\boldsymbol\ast'$, we have $F(a)\not\subseteq\bigcup_{i\in\gamma+1}B_i$. Hence, for every $\gamma\in{\sf Ord}$, there exists some $c_\gamma\in F(a)$ such that $c_\gamma\not\in\bigcup_{i\in\gamma+1}B_i\;(\ddag)$. By $F(a)=\bigcup_{\alpha\in\textsf{Ord}}B_\alpha$, there is some $\delta\in{\sf Ord}$ such that $c_\gamma\in B_\delta$, and so by $\ddag$, we have $\gamma\in\delta$ (since $\delta\not\in\gamma+1$). Therefore, for each ordinal $\gamma$, we have  $c_\gamma\in B_\delta\setminus\bigcup_{i\in\gamma+1}B_i$ for some bigger ordinal $\delta$($\ni\gamma$). Thus, the set $\mathcal{B}$ exhausts the class ${\sf Ord}$, which is impossible.
\end{proof}

We end the paper with two observations about the set $\mathcal{B}$ above. First, it can be defined as $\mathcal{B}=\bigcup_{i\in I}B_i$, where $(I,\prec,B)$ is a $\sqsubseteq$-maximal set such that for some well-ordering $\prec$ on $I$ and some assignment $i\mapsto B_i\in\mathscr{P}(A)$, the following holds:  $$(\boldsymbol\ast'')\;\forall j\in I: B_j=\{x\in A\mid F(x)\subseteq\bigcup_{i\prec j}B_i\}.$$
Second, we prove in the appendix that it is equal to the set $$\mathcal{D}_{\infty}=\{a\in A\mid \boldsymbol\not\!\exists
\{x_i\}_{i=1}^{\infty}\subseteq A\!\!: x_1\in F(a)\wedge\bigwedge_{n\in\mathbb{N}}x_{n+1}\in F(x_n)\},$$
which was directly shown in \textbf{\cite{R05}} to be out of the $F$'s range.

\appendix
\section*{\textsc{Appendix}}

\begin{proposition}[$\mathcal{B}=\mathcal{D}_\infty$]
\,

\noindent
The set $\mathcal{B}$, defined by $(\boldsymbol\ast')$ above, is equal to the set $\mathcal{D}_\infty$ (defined above also).
\end{proposition}

\begin{proof}

\noindent
$\boldsymbol( \mathcal{B} \subseteq \mathcal{D}_\infty \boldsymbol)$: If $x\in\mathcal{B}$ but $x\not\in\mathcal{D}_\infty$, then for a sequence $\{a_i\}_{i=1}^{\infty}\subseteq A$, we have $a_1\in F(x)$ and $a_{j+1}\in F(a_j)$ for each $j=1,2,\cdots$. By $x\in\mathcal{B}$, we have $x\in B_{\alpha_0}$ for some $\alpha_0\in{\sf Ord}$, and so by $\boldsymbol\ast'$, we have $F(x)\subseteq\bigcup_{i\in\alpha_0}B_i$. Thus, $a_1\in B_{\alpha_1}$ for some $\alpha_1\in\alpha_0$. Continuing this way, by $\boldsymbol\ast'$, we have
$F(a_1)\subseteq\bigcup_{i\in\alpha_1}B_i$,
thus $a_2\in B_{\alpha_2}$ for some $\alpha_2\in\alpha_1$, and so on and so forth. At the end, we get an infinite decreasing sequence of ordinals $\cdots\in\alpha_2\in\alpha_1\in\alpha_0$, which is a contradiction.

\noindent
$\boldsymbol( \mathcal{D}_\infty \subseteq \mathcal{B} \boldsymbol)$: If $x\not\in\mathcal{B}$, then $F(x)\not\subseteq\mathcal{B}$, since otherwise if for each $y\in F(x)$ there existed some $\alpha_y\in{\sf Ord}$ such that $y\in B_{\alpha_y}$, then for an ordinal $\gamma$ which is bigger than all the ordinals in the set $\{\alpha_y\mid y\in F(x)\}$, we would have $F(x)\subseteq\bigcup_{i\in\gamma}B_i$, thus $x\in B_\gamma$ by $\boldsymbol\ast'$, so $x\in\mathcal{B}$, contradicting the assumption (of $x\not\in\mathcal{B}$). So, for some $x_1\in F(x)$, we have $x_1\not\in\mathcal{B}$. Thus, $F(x_1)\not\subseteq\mathcal{B}$ again, and so there should exist some $x_2\in F(x_1)$ such that $x_2\not\in\mathcal{B}$. Repeating this will produce a sequence $\{x_i\}_{i=1}^\infty\subseteq A$ with $x_1\in F(x)$ and $x_{j+1}\in F(x_j)$ for each $j=1,2,\cdots$. Thus, $x\not\in\mathcal{D}_\infty$.
\end{proof}

\end{document}